\documentclass[12pt,a4paper,twoside]{article}
\usepackage[english]{babel}
\usepackage{amssymb}
\usepackage{inputenc}
\usepackage{amssymb}
\begin{document}

\title{\bf Innerness of Derivations on Subalgebras of Measurable Operators}

\author{Sh. A. Ayupov  $^{1 *},$  K. K. Kudaybergenov  $^2$}

\maketitle

\begin{abstract}

Given a  von Neumann algebra  $M$ with a faithful normal
semi-finite trace $\tau,$ let $L(M, \tau)$ be the algebra of all
$\tau$-measurable operators affiliated with $M.$  We prove that if
$A$ is a locally convex reflexive complete metrizable solid
$\ast$-subalgebra in $L(M, \tau),$ which can be embedded into a locally
bounded weak Fr\'{e}chet $M$-bimodule, then any derivation on  $A$ is
inner.

\end{abstract}

\medskip
$^1$ Institute of Mathematics and information  technologies,
Uzbekistan Academy of Science, F. Hodjaev str. 29, 100125, Tashkent
(Uzbekistan), e-mail: \emph{sh\_ayupov@mail.ru,
e\_ayupov@hotmail.com}

 $^{2}$ Institute of
Mathematics and information  technologies, Uzbekistan Academy of
Science, F. Hodjaev str. 29, 100125, Tashkent (Uzbekistan), e-mail:
\emph{karim2006@mail.ru}

\medskip \textbf{AMS Subject Classifications (2000): 46L57, 46L50, 46L55,
46L60}
\\

\textbf{Key words:}  von Neumann algebras,    measurable operator,  weak Fr\'{e}chet  bimodule,
 derivation, inner derivation.

* Corresponding author
\newpage

\section*{\center 1. Introduction}

 The structure of automorphisms and derivations of operator algebras
 is an important part of the theory of operator algebras and their
  applications in quantum dynamics.

   Recall that a linear operator $D$ on an algebra
    $A$ is called a  \emph{derivation} if it satisfies the condition
    $$D(xy)=D(x)y+xD(y)$$ for all $x, y\in A.$

    Every  (but fixed) element $a\in A$ generates a derivation $D_a$ on $A,$
   defined as $D_a(x)=ax-xa,\,  x\in A.$ Such derivations are said to be  \emph{inner} derivations.

   Derivations on $C^{\ast}$-algebras and von Neumann  algebras have been studied in the monographs
   of Sakai  \cite{Sak1},  \cite{Sak2}.
  It is well-known that every derivation on a
  $C^{\ast}$-algebra $A$ is norm continuous and if $C^{\ast}$-algebras $A$ is unital and
  simple or it is weakly closed (i.e.
  is  a von
 Neumann algebra) then any derivation on $A$ is  inner.
 For general
 Banach algebras similar problems were considered in the  monograph \cite{Dal}.

Investigation of derivations on unbounded operator algebras and, in particular,
on the algebra $L(M)$ of measurable operators affiliated with a
von Neumann algebra $M,$ was initiated in the papers  \cite{Ayu1}, \cite{Ayu2}.
  One of the main problems posed in these papers was: whether any derivation
  on $L(M)$ is inner.
   A negative answer to this problem in the general setting was given in the paper
    \cite{Ber} (see also \cite{Kus}).  Namely, it was proved that if $M$
    is a non atomic abelian von Neumann algebra
    (in particular   $L^{\infty}(0; 1)$) then $L(M)$ (resp.  $L^{0}(0; 1)$) admits a non trivial
     (and hence discontinuous, and non inner) derivation.

     Further there were some positive results on this way. In
       \cite{Alb1} we have proved that if $M$ is  type I  von Neumann algebra with a faithful normal
semi-finite trace $\tau,$ then a derivation on the algebra
  $L(M, \tau)$ of all $\tau$-measurable operators affiliated with $M$ is inner if and only if it is $Z$-linear,
or equivalently if it is identically zero on the center $Z$ of $M.$
Recently \cite{Alb3} we gave a complete description of derivation on $L(M,
\tau)$ and in particular proved that if $M$ is of type I$_\infty$ then any
derivation on $L(M, \tau)$ is inner.

 If $M$ is a general von Neumann algebra with a faithful normal
semi-finite trace $\tau,$ then the algebra $L(M, \tau)$ contains various subalgebras
 with different properties of derivations.
One of the interesting classes of  subalgebras in $L(M, \tau)$ are so called
 Arens algebras
      $$L^{\omega}(M, \tau)=\bigcap\limits_{p\geq1}L^{p}(M, \tau).$$
        In  the paper \cite{Alb} we gave a  complete description of derivations on $L^{\omega}(M, \tau)$ and proved that any derivation on $L^{\omega}(M, \tau)$
        is inner if and only if the trace  $\tau$ is finite.

        In this connection a natural question arises:

        which subalgebras in $L(M, \tau)$ admit only inner derivation?

In this  paper we give a sufficient condition for subalgebras in  $L(M, \tau)$ to
have such a property. Namely, we prove that if   $A$ is locally
convex reflexive complete metrizable solid $\ast$-subalgebra in
$L(M, \tau),$ which can be embedded into a locally bounded weak Fr\'{e}chet
$M$-bimodule, then any derivation on  $A$ is inner.

\begin{center}  \textbf{2. Preliminaries}
\end{center}

Let $H$ be a Hilbert space, and let  $B(H)$ be the algebra of all
bounded linear operators on   $H.$ Consider a von Neumann  algebra
$M\subset B(H)$ with a faithful normal semi-finite trace $\tau,$ and
denote by $\mathcal{P}(M)$ the lattice of (orthogonal) projections
in $M.$

A linear subspace  $\mathcal{D}$ in  $H$ is affiliated with  $M$ (denoted as
$\mathcal{D}\eta M$), if $u(\mathcal{D})\subset \mathcal{D}$ for any unitary operator $u$ from the
commutant $$M'=\{y\in B(H):xy=yx, \,\forall x\in M\}$$ of the von
Neumann algebra $M.$

A linear operator  $x$  with the domain  $\mathcal{D}(x)\subset H$ is said to
be affiliated with  $M$ (denoted as  $x\eta M$) if $u(\mathcal{D}(x))\subset
\mathcal{D}(x)$ and $ux(\xi)=xu(\xi)$ for all $u\in M',$ $\xi\in \mathcal{D}(x).$

A linear subspace  $\mathcal{D}$ in $H$ is called  $\tau$-dense, if

1) $\mathcal{D}\eta M;$

2) given any  $\varepsilon>0$ there exists a projection
$p\in\mathcal{P}(M)$ such that  $p(H)\subset \mathcal{D}$ and
$\tau(p^{\perp})\leq\varepsilon.$

A closed linear operator  $x$ is called $\tau$-\emph{measurable}
with respect to the von Neumann algebra $M,$ if $x\eta M$ and $\mathcal{D}(x)$
is  $\tau$-dense in $H.$

 Denote by $L(M,
\tau)$ the set of all $\tau$-measurable operators affiliated with
 $M.$ Consider the topology  $t_{\tau}$ of convergence in measure on $L(M, \tau),$ which is defined by
 the following neighborhoods of zero:
$$V(\varepsilon, \delta)=\{x\in L(M, \tau): \exists e\in\mathcal{P}(M), \tau(e^{\perp})\leq\delta, xe\in
M,  \|xe\|_{M}\leq\varepsilon\},$$ where $\|\cdot\|_{\infty}$ is the
operator norm on $M,$ and  $\varepsilon, \delta$ are positive numbers.

 It is well-known
\cite{Nel} that $L(M, \tau)$ equipped with the measure topology is a
complete metrizable topological $\ast$-algebra.

Now let us recall the notion of a bimodule over a Banach algebra
(see \cite{Dal}).

Let  $A$ be a complex algebra and  let $E$ be a complex
linear space.  $E$ is called a left $A$-module
(respectively right $A$-module) if a bilinear map $(a,
x)\mapsto a\cdot x$ (respectively $(a, x)\mapsto x\cdot a$) from
$A\times E$ into $E$ is defined,
 such that  given any  $a, b\in A$ and $x\in E$ one has
  $$a\cdot(b\cdot
x)=ab\cdot x \quad (\mbox{respectively}\, (x\cdot a)\cdot b=x\cdot
ab),$$

$E$ is said to be   $A$-bimodule if   is  left and right
 $A$-module simultaneously, and
 $$a\cdot(x\cdot b)=(a\cdot x)\cdot b,$$ for all  $a, b\in
A, x\in E.$

 Let  $A$ be a Banach algebra and suppose that $E$ is a
  Fr\'{e}chet  space, i.e. a complete metric space
 with a shift invariant metric. If $E$ is an $A$-bimodule and the maps  $x\mapsto
a\cdot x$ and $x\mapsto x\cdot a$ are continuous for each $a\in A,$
then $E$ is called a \emph{weak Fr\'{e}chet}  $A$-bimodule.

Interesting examples of weak Fr\'{e}chet $A$-bimodules are given by non
commutative $L^{p}$-spaces  $L^{p}(M, \tau)\subset L(M, \tau),$
$p\geq1.$ Indeed, given any  $a\in M$ and $x\in L^{p}(M, \tau)$ one
has $ax\in L^{p}(M, \tau),$ $xa\in L^{p}(M, \tau)$ and $ \|a
x\|_{p}\leq \|a\|_{\infty}\|x\|_p$ which imply the above statement.

Let $E$ and $F$ be metrizable linear topological spaces and let
$T:E\rightarrow F$ be a linear operator. The \emph{separating space}
of the linear map $T,$ denoted by $\mathcal{S}(T),$ is defined as
$$\mathcal{S}(T)=\{y\in F:\mbox{there is}\,\,(x_{n})_{n\in\mathbb{N}}
 \,\,\mbox{in} \,\, E
 \,\,\mbox{such that}\,\, x_{n}\rightarrow 0\,\,\mbox{and}\,\, T(x_{n})\rightarrow y\}.$$

Recall that $\mathcal{S}(T)$ is closed (see \cite{Dal},
{Proposition 5.1.2}) and that the closed graph theorem is valid
for complete metrizable topological linear space. Therefore $T$ is
continuous if and only if $\mathcal{S}(T)=\{0\}.$

 \textbf{Definition 2.1.} \cite{Dal}.
Let $A$ be an algebra, and let $E$ be a topological linear space
which is a $A$-bimodule. Then $E$ is a \emph{separating module} if, for
each sequence $\{a_n\}$ in $A,$ the nest $\overline{(a_1\cdots a_n
E)}$ stabilizes, i.e. there is a $n_0\in\mathbb{N}$ such that
$\overline{(a_1\cdots a_n E)}=\overline{(a_1\cdots a_{n+1} E)}$
for all $n>n_0,$ where $\overline{F}$ -- the  closure of the set
$F.$

 A linear topological space  $E$ is said to be \emph{locally bounded} if there exists a bounded
 neighborhood of zero in $E.$

From \cite[Theorem 5.2.15]{Dal}, we have the following

 \textbf{Proposition 2.2.} \emph{Let  $E$ be a
weak Fr\'{e}chet $A$-bimodule, and let $D:A\rightarrow E$ be a
derivation. Then}

1) \emph{The separating space $\mathcal{S}(D)$ is a closed submodule of
$E;$}

2) \emph{Suppose that $E$ is locally bounded. Then $\mathcal{S}(D)$ is a
separating module.}

Given a linear topological space  $X$ with a topology  $t_X,$ let us
denote by  $x_n\stackrel{t_X}{\longrightarrow}0$ the convergence in
the topology $t_X.$

If  $(A, t_A)$ and  $(B, t_B)$ are linear topological spaces, with
$A\subseteq B\subseteq L(M, \tau)$ then we shall suppose that the
topology  $t_A$ is stronger than $t_B,$ i.e.  $(A, t_A)$ is
topologically imbedded into $(B, t_B).$

\begin{center}
{\bf 3. The main results}
\end{center}

 The aim of the present section is to prove the following result.

\textbf{Theorem 3.1} \emph{Let  $M$  be a von Neumann algebra with a
faithful normal semi-finite trace $\tau.$ Suppose that $A$ is a
complete metrizable solid $\ast$-subalgebra in  $L(M, \tau)$ and $E$
is a locally bounded  weak Fr\'{e}chet $M$-bimodule in $L(M, \tau).$
If}

1) \emph{$A$ is locally convex and reflexive;}

2) \emph{$M\subset A\subset E$ are topological imbedding,}

\emph{then any derivation of the algebra $A$ is inner.}

Recall that   a subalgebra $A$ in   $L(M, \tau)$ is solid if
$x\in A,$ $y\in L(M, \tau), |y|\leq |x|$ implies $y\in A.$

The proof of this theorem consists of several steps.

\textbf{Proposition 3.2.} \emph{Given  an arbitrary  von Neumann
algebra $M,$ and  a weak Fr\'{e}chet $M$-bimodule $E,$ suppose that
$p$ is a projection in $M$ and  $D:M\rightarrow E$ is a derivation,
i.e. a linear map such that $D(xy)=D(x)y+xD(y)$ for all $x, y\in M.$
Put $D_p(x)=pD(x)p,\,x\in pMp.$
 Then}

1) $D_p:pMp\rightarrow pEp$ \emph{is a derivation;}

2) $p\mathcal{S}(D)p\subseteq\mathcal{S}(D_p).$

Proof. 1) For   $x, y\in pMp$ we have $x=pxp, y=pyp.$ Therefore
$D_p(xy)=pD(pxyp)p=pD(pxppyp)p=pD(pxp)pyp+pxpD(pyp)p=D_p(x)y+xD_p(y),$
i.e. $D_p$ is a derivation.

2) For  $y\in\mathcal{S}(D)$ according the definition there exists a
sequence  $\{x_n\}$ in $M$ such that
$x_n\stackrel{\|\cdot\|}{\longrightarrow}0$ and
$D(x_n)\stackrel{t_E}{\longrightarrow} y$ as $n\rightarrow\infty.$
But then $px_np\stackrel{\|\cdot\|}{\longrightarrow}0$ and
$D_p(px_np)=pD(px_n
p)p=pD(p)x_np+pD(x_n)p+px_nD(p)p\stackrel{t_E}{\longrightarrow}
pyp,$ i.e. $D_p(px_np)\stackrel{t_E}{\longrightarrow}pyp,$ which
means that $pyp\in \mathcal{S}(D_p)$ and therefore
$p\mathcal{S}(D)p\subseteq \mathcal{S}(D_p).$ The proof is complete.
$\blacksquare$

 \textbf{Proposition 3.3.} \emph{Let $M$ be a von Neumann algebra with a
faithful normal semi-finite trace $\tau$ and let $E$ be a locally
bounded  weak Fr\'{e}chet $M$-bimodule in $L(M, \tau).$ Then every
derivation  $D:M\rightarrow E$ is automatically continuous.}

Proof. Let us show that $\mathcal{S}(D)=\{0\}.$ Suppose the
opposite, i.e. $\mathcal{S}(D)\neq\{0\}$ and take a non zero  $y\in
\mathcal{S}(D).$ Chose a projection  $e$ in $M$ such that $ye\in M$
and $ye\neq0.$ Since  $\mathcal{S}(D)$ is a submodule in  $E,$ we have
that  $ye(ye)^{\ast}\in \mathcal{S}(D).$ Thus without loss of the
generality we may suppose that $y\geq0.$ Take a projection  $p\in
M$ such that $pyp\neq0$ and $n^{-1}p\leq pyp \leq np$ for an
appropriate  $n\in\mathbb{N}.$ Then $pyp$ is invertible in $pMp,$
i.e. there exists an element  $z\in pMp$ such that $pypz=p.$ Since
$\mathcal{S}(D)$ is an $M$-bimodule it follows that  $p\in
\mathcal{S}(D).$ Now consider two cases separately:

\emph{The case 1.} $pMp$ is finite dimensional. Observe the
derivation  $D_p:pAp\rightarrow pEp$ defined by
$$D_p(x)=pD(pxp)p,\quad x\in pAp.$$
Since  $pMp$ is finite dimensional, the spaces $pAp$ and $pEp$ are
also finite dimensional as subspaces of $L(M, \tau)=pMp.$ Therefore
 $D_p$ is necessary continuous, i.e.  $\mathcal{S}(D_p)=\{0\}.$

On the other hand Proposition 3.2 implies that
$p\mathcal{S}(D)p\subset \mathcal{S}(D_p),$ and from the
consideration above we have that $p\in \mathcal{S}(D_p)$ and hence
$0\neq p\in p\mathcal{S}(D)p\subset\mathcal{S}(D_p)=\{0\}.$ This
contradiction implies that  $\mathcal{S}(D)=\{0\}.$

\emph{The case 2.} $pMp$ is infinite dimensional. In this case there
exists a strictly monotone decreasing sequence  $(p_n)$ of
projections in $M$ such that  $p_n\leq p$ for all $n\in\mathbb{N}.$
Since  $\mathcal{S}(D)$ is a closed submodule in  $E$ (Proposition
2.2)  $p_n\mathcal{S}(D)$ is also closed. Further  $p_n\neq
p_{n+1}\leq p$ implies that
$$\overline{p_n\mathcal{S}(D)}\neq\overline{p_{n+1}\mathcal{S}(D)}.\eqno (1)$$
On the other hand, since  $E$ is locally bounded Proposition 2.2
implies that the nest  $\{\overline{p_n\mathcal{S}(D)}\}$ stabilizes,
i.e. there exist  $n_0\in\mathbb{N}$ such that
$\overline{p_n\mathcal{S}(D)}=\overline{p_{n+1}\mathcal{S}(D)}$ for
all $n>n_0$ in a contradiction with  (1). Therefore
$\mathcal{S}(D)=\{0\}$ and the proof is complete.

\textbf{Remark 1.} Proposition 3.3  implies the well-known fact that
every derivation on the von Neumann algebra $M$ is norm continuous.
It is sufficient to put  $E=M.$

\textbf{Proposition 3.4.} \emph{Let $M$ be a von Neumann algebra
with a faithful normal semi-finite trace $\tau.$ Suppose that  $E$
is a locally bounded  weak Fr\'{e}chet $M$-bimodule in $L(M, \tau)$
and  $A$ is complete metrizable algebra such that  $A\subseteq E.$
Then any derivation  $D:M\rightarrow A$ is continuous.}

Proof. Let $\{x_n\}\subset M,
x_n\stackrel{\|\cdot\|}{\longrightarrow}0$ and
$D(x_n)\stackrel{t_A}{\longrightarrow}y,$ which implies that
$D(x_n)\stackrel{t_\tau}{\longrightarrow} y.$ Let us show that
$y=0.$ By Proposition 3.3 $D:A\rightarrow E$ is continuous and thus
$D(x_n)\stackrel{t_E}{\longrightarrow} 0$ and hence
$D(x_n)\stackrel{t_\tau}{\longrightarrow}
 0.$ This implies that $y=0.$ The
proof is complete. $\blacksquare$

If  $A$ is a locally convex metrizable space, then  its topology can
be generated by an increasing sequence of seminorms $\{\rho_n,
n\in\mathbb{N}\}.$

Denote by $U$ the group of all unitaries of the von Neumann algebra
 $M.$

\textbf{Proposition  3.5.} \emph{Let $A$ be a locally convex weak
Fr\'{e}chet $M$-bimodule in $L(M, \tau).$ Given any non zero element
$x\in A$ there exists a seminorm $\rho_n$ such that}
$$\inf\{\rho_n(uxu^{\ast}):u\in U\}\neq0.$$

Proof. Suppose that opposite, i.e.  $\inf\{\rho_n(uxu^{\ast}):u\in
U\}=0$ for all $n\in\mathbb{N}.$ Chose the unitaries  $u_n\in U$
such that  $\rho_n(u_nxu_n^{\ast})\leq n^{-1}.$ Since
$\rho_k\leq\rho_{k+1},$ we have that
$\rho_k(u_nxu_n^{\ast})\rightarrow0$ as $n\rightarrow\infty$ for
each fixed  $k\in\mathbb{N}.$ This means that
$u_nxu_n^{\ast}\stackrel{t_A}{\longrightarrow}0,$ É and hence
$u_nxu_n^{\ast}\stackrel{t_\tau}{\longrightarrow}0.$ Since
$\|u_n\|_{\infty}=1$ for all $n=1,2,...,$ we obtain that
$x=u_n^{\ast}(u_nxu_n^{\ast})u_n\stackrel{t_\tau}{\longrightarrow}0,$
i.e.  $x=0$ a contradiction.  The proof is complete. $\blacksquare$

\textbf{Proposition 3.6.} \emph{Let  $M,$ $A$ and $E$ be as in
theorem 3.1. Then every derivation  $D:M\rightarrow A$ is spatial, i.e. $D(x)=ax-xa$
for an appropriate $a\in A$ and every $x\in M$.}

Proof.  By Proposition 3.4 the derivation $D:M\rightarrow A$ is
continuous.

Let $U$ be the group of all unitary elements in $M.$ Given any   $u\in U$
put
$$T_u(x)=uxu^{\ast}+D(u)u^{\ast}, \quad x\in A.$$

Since map $x\mapsto uxu^{\ast}$ is continuous,  the map  $T_u$  is
$\sigma(A, A^{\ast})$-continuous.

For $u, v\in U$ we have
\begin{center}
$T_u(T_v(x))=T_u(vxv^{\ast}+D(v)v^{\ast})=u(vxv^{\ast}+D(v)v^{\ast})u^{\ast}+D(u)u^{\ast}=
uvxv^{\ast}u^{\ast}+uD(v)v^{\ast}u^{\ast}+D(u)u^{\ast}=(uvx+D(u)v+uD(v))(uv)^{\ast}=
uvx(uv)^{\ast}+D(uv)(uv)^{\ast}=T_{uv}(x),$
\end{center}
i. e.  $$ T_uT_v=T_{uv}, \quad u, v\in U.\eqno (2)$$

Further we have  $T_u(0)=D(u)u^{\ast}$ and the  continuity of $D$ implies that the  set
$K_D=\{T_v(0): v\in U\}=\{D(v)v^{\ast}: v\in U\}$ is bounded in $A.$ Moreover, the set
$K=\mbox{cl(co}(K_d))$ -- the  closure of the convex hull of
$K_D,$ is a closed convex bounded subset in  $A.$ The reflexivity
of the space  $A$ then implies that  $K$ is a non-void  $\sigma(A, A^{\ast})$-compact convex
set. From (2) it follows that  $T_u(K_D)\subset K_D$ for all $u\in U.$ Since
$T_u$ is an affine homeomorphism we have
$T_u(\mbox{cl(co}(K_D)))=\mbox{cl(co}(T_u(K_D)))\subset
\mbox{cl(co}(K_D)),$ i. e. $T_u (K)\subseteq K$
for all $u\in U.$

According  Proposition  3.6, given any  $x, y\in A, x\neq y$ there exists seminorm
$\rho_n$ such that
$$\inf\{\rho_n(T_u(x)-T_u(y)):u\in U\}=\inf\{\rho_n(u(x-y)u^{\ast}):u\in U\}\neq0.$$

Therefore  $\{T_u: u\in U\}$ is a non-contracting (in the sense of \cite{Nam}) semigroup
of $\sigma(A, A^{\ast})$-continuous affine mappings of a $\sigma(A, A^{\ast})$-compact convex
set $K.$ By
Ryll-Nardzewski's fixed point theorem \cite{Nam}, there exists $a\in
K$ such that  $T_u (a)=a$ for all $u\in U.$ This means that
$uau^{\ast}+D(u)u^{\ast}=a,$ i. e. $D(u)=au-ua$ for all $u\in U.$
Since every element of $M$ is a  linear combination of unitaries
from $M$, we have $D(x)=ax-xa$ for all $x\in M,$ i.e. $D$ is a spatial derivation on $M$
with values in $A.$ The proof is
complete. $\blacksquare$

Proof of Theorem 3.1.  According to Proposition 3.6 there exists
element $a\in A$ such that  $D(x)=ax-xa$ for all $x\in M.$ We shall proof that this is true
 for all $x\in A.$

 First suppose that $x\in A,\, x\geq0.$
In this case  the element $\textbf{1}+x\in A\subset L(M, \tau)$ is  invertible and
moreover $(\textbf{1}+x)^{-1}\in M.$

For an invertible element
   $b\geq0$ in $A$ one has
$0=D(\textbf{1})=D(bb^{-1})=D(b)b^{-1}+bD(b^{-1}),$ i. e.
$D(b)=-bD(b^{-1})b.$

Therefore
$$D(x)=D(\textbf{1}+x)=-(\textbf{1}+x)D((\textbf{1}+x)^{-1})(\textbf{1}+x).$$
 On  the other hand, since   $(\textbf{1}+x)^{-1}\in M$ we have from the above
$$D((\textbf{1}+x)^{-1})=a(\textbf{1}+x)^{-1}-(\textbf{1}+x)^{-1}a.$$
Therefore,
$$-(\textbf{1}+x)D((\textbf{1}+x)^{-1})(\textbf{1}+x)=
-(\textbf{1}+x)[a(\textbf{1}+x)^{-1}-$$$$-(\textbf{1}+x)^{-1}a](\textbf{1}+x)=-(\textbf{1}+x)a+a(\textbf{1}+x)=
ax-xa,$$  i. e.
$$D(x)=-(\textbf{1}+x)D((\textbf{1}+x)^{-1})(\textbf{1}+x)=ax-xa.$$
Therefore $D(x)=ax-xa$ for every $x\in A,\,x\geq0.$
Since   $A$ is a solid $\ast$-subalgebra in $L(M, \tau),$  every
element from   $A$  is a linear combination of positive
elements of $A.$ Thus  $D(x)=ax-xa$ for all $x\in A.$ The
proof is complete. $\blacksquare$

 \textbf{Remark  2.}  The condition on  $A$ to be solid is used only at the end of proof of Theorem
 3.1. In fact this condition may be replaced by the condition
  $Lin(A_+)=A,$ where $A_+$ is the positive cone of  $A$ and
  $Lin(A_+)$ is the linear span of $A_+$ (i.e. the positive cone $A_+$ is hereditary in
  $A$).

\textbf{Example 3.8.} An example of algebras, satisfying the
conditions of Theorem  3.1 is given by a non commutative Arens
algebra  $L^{\omega}(M,\tau)$ in the case of a finite trace $\tau$
(see \cite{Alb}).

Given  be a von Neumann algebra $M$ with a
faithful normal semi-finite trace $\tau$ $M$ and  $p\geq1,$ put
$L^{p}(M, \tau)=\{x\in L(M, \tau):\tau(|x|^{p})<\infty\}.$ It is
known  \cite{Nel} that  $L^{p}(M, \tau)$ is a Banach space with
respect to the norm
$$\|x\|_p=(\tau(|x|^{p}))^{1/p},\quad x\in L^{p}(M, \tau).$$
Consider the space
\begin{center}
$L^{\omega}(M, \tau)=\bigcap\limits_{p\geq1}L^{p}(M, \tau).$
\end{center}
It is known \cite{Abd}, \cite{Ino}  that $L^{\omega}(M, \tau)$ is a
locally convex metrizable $\ast$-algebra with the topology $t$
generated by the sequence of norms
$$\|x\|_n'=\max\{\|x\|_1, \|x\|_n\},\,n\in\mathbb{N}.$$
The algebra $L^{\omega}(M, \tau)$ is called a  (non commutative)
\emph{Arens algebra}. The dual space for $(L^{\omega}(M, \tau), t)$
was described in  \cite{Abd}, where it has been proved that
$(L^{\omega}(M, \tau), t)$ is reflexive if and only if trace $\tau$
is finite.

Therefore Theorem  3.1 implies that if the trace $\tau$ is finite,
then every derivation on the algebra  $L^{\omega}(M, \tau)$ is
inner.

It should be noted also that a complete description of derivations
on general $L^{\omega}(M, \tau)$ was obtain in \cite{Alb}. Namely.
it has been proved that every derivation on $L^{\omega}(M, \tau)$ is
spatial and has the form
$$D(x)=ax-xa, \quad x\in L^{\omega}(M, \tau)$$
  for an appropriate  $a\in M+L^{\omega}_2(M, \tau),$ where $L^{\omega}_{2}(M,
\tau)=\bigcap_{p\geq2}L^{p}(M, \tau).$

Now let us consider an example of an algebra $A$ satisfying all
conditions of Theorem 3.1 except $M\subset A,$ which admits non-inner
derivations.

 \textbf{Example 3.9.} Put
 $$A=L^{\omega}_{2}(M, \tau)=\bigcap_{p\geq2}L^{p}(M, \tau)$$ and consider
 on $A$ the topology generated by the system of norms  $\{\|\cdot\|_{p\geq2}\}.$
 Then $A$ is a metrizable locally convex  $\ast$-algebra \cite{Alb} and
$$L^{\omega}_{2}(M,
\tau)=\bigcap_{n=2}^{\infty}L^{n}(M, \tau).$$ As the intersection of
countable family of reflexive Banach space,  $A$ is also
reflexive. If the trace $\tau$ is semi-finite but not finite, $M$ is
not contained in  $A=L^{\omega}_2(M, \tau).$ Every derivation of
$L^{\omega}_{2}(M, \tau)$ has the form
$$D_a(x)=ax-xa,\, x\in L^{\omega}_{2}(M, \tau)$$
for an appropriate  $a\in M+L^{\omega}_{2}(M, \tau)$ (see \cite{Alb}).

Now  if $M$ is non commutative and $a$ is a non central element from
$M\setminus L^{\omega}_2(M, \tau)$ the spatial derivation  $D_a$ on
$L^{\omega}_{2}(M, \tau)$ is not inner.

Let us consider an example of an algebra $A$ which shows that the
reflexivity of $A$ is not a necessary condition for the statement of
Theorem 3.1.

 \textbf{Example 3.10.} Let $M$ be the  $C^{\ast}$-product of von Neumann algebra
  $M_n,$ i.e.
$$M=\bigoplus\limits_{n=1}^{\infty}M_n=\{\{x_n\}: x_n\in M_n,\,\sup\|x\|_{M_n}<\infty\},$$
where $\|\cdot\|_{M_n}$ is the $C^{\ast}$-norm on $M_n.$

Put
$$A=\{\{x_n\}: x_n\in M_n,\,\sum\limits_{n=1}^{\infty}\frac{\textstyle 1}{\textstyle 2^{n}}\|x_n\|_{M_n}
^{p}<\infty,\,1\leq p<\infty\}.$$ Then it is clear that $\{x_n\}\in
A$ if and only if  $\{\|x_n\|_{M_n}\}\in l^{\omega},$ where
$l^{\omega}=\bigcap\limits_{p\geq1}l^{p}$ -- the Arens algebra
associated with abelian von Neumann algebra  $l^{\infty}$ of all
bounded complex sequences with the trace
$$\nu(\{\lambda_n\})=\sum\limits_{n=1}^{\infty}\frac{\textstyle
1}{\textstyle 2^{n}}\lambda_n, \{\lambda_n\}\in l^{\infty}.$$

 Consider the topology  $t$ on  $A,$ generated by the family of norms
 $$\|x\|_{A, p}=(\sum\limits_{n=1}^{\infty}\frac{\textstyle 1}{\textstyle 2^{n}}\|x_n\|_{M_n}
^{p}) ^{\frac{1}{p}}.$$ With the coordinatewise algebraic operations
and involution $(A, t)$ becomes a locally convex complete metrizable
 $\ast$-algebra. If at least one of the algebras $M_n$ is infinite dimensional,
  the  $A$ is not isomorphic to any Arens algebra.

Consider a derivation $D:A\rightarrow A.$ Let  $q_n$ be central
projection in $M$ such that   $q_nM=M_n.$ Then we have
$$D(q_n x)=q_n D(x),\quad x\in A,$$
and therefore  $D(M_n)\subseteq M_n$  and the restriction
$$D_n(x)=q_n D(x),\quad x\in M_n,$$
 gives a derivation  $D_n:M_n\rightarrow M_n.$ The classical theorem
 of Sakai   implies the existence of an appropriate
  $a_n\in M_n$ such that $D_n(x)=a_n x- x a_n $ for all
$x\in M_n,$ moreover one can assume that $\|a_n\|\leq \|D_n\|$ (see \cite[Theorem 4.1.6]{Sak1})

Let us show that  $a=\{a_n\}\in A$ and that $D(x)=ax- xa$ for all
$x\in A.$

Take an arbitrary  $\varepsilon>0.$  For each $n\in\mathbb{N}$ there
exists  $x_n\in M_n,\,\|x_n\|_{M_n}\leq1$ such that
$\|D_n\|\leq\|D_n(x_n)\|+\varepsilon.$ Then  $x=\{x_n\}\in M$ and
$D(x)=\{D_n(x_n)\}\in A.$ Therefore $\{\|D(x_n)\|_{M_n}\}\in
l^{\omega}.$ á The inequalities
$$\sum\limits_{n=1}^{\infty}\frac{\textstyle 1}{\textstyle 2^{n}}\|a_n\|_{M_n}
^{p}\leq\sum\limits_{n=1}^{\infty}\frac{\textstyle 1}{\textstyle
2^{n}}\|D_n\|^{p}\leq\sum\limits_{n=1}^{\infty}\frac{\textstyle
1}{\textstyle
2^{n}}(\|D_n(x_n)\|+\varepsilon)^{p}\leq$$$$\leq2^{p-1}\sum\limits_{n=1}^{\infty}\frac{\textstyle
1}{\textstyle 2^{n}}(\|D_n(x_n)\|^{p}+\varepsilon^{p})$$ imply that
$\{\|a_n\|_{M_n}\}\in l^{\omega},$ i.e. $a\in A.$ Further since
$$q_n D(x)=D_n(q_n x)=a_n q_n x - q_n x a_n=q_n(ax-xa)$$ for all $n=1,2,...,$ taking the sum over all $q_n$
 we obtain  $D(x)=ax-xa$ for all $x\in A.$ The proof is
complete. $\blacksquare$

 It well-known that every abelian von Neumann algebra is isomorphic to an algebra
  $L^{\infty}(\Omega)$ of all essentially bounded measurable complex functions on a measure space
  $(\Omega, \Sigma, \mu).$  In this case the algebra $\tau$-measurable operators
  $L(M, \tau)$
 is isomorphic with the algebra  $L^{0}(\Omega)$ of all measurable functions on $\Omega.$

\textbf{Proposition  3.11.} \emph{Let  $A$ be a $\ast$-subalgebra in
$L^{0}(\Omega)$ and let $E$ be a locally bounded weak Fr\'{e}chet
$L^{\infty}(\Omega)$-bimodule such that $L^{\infty}(\Omega)\subseteq
A\subseteq E.$ Then every derivation on  $A$ is identically zero.}

Proof. Since  $L^{\infty}(\Omega)$ is abelian, every derivation  $D$
is equal to zero on idempotents (projections) from
$L^{\infty}(\Omega).$  Therefore  $D(x)=0$ on each step function
$x\in A.$ The space of step functions is dense in
$L^{\infty}(\Omega)$ and by Proposition  3.4 the derivation
$D:L^{\infty}(\Omega)\rightarrow A$ is continuous, therefore
$D(x)=0$ for all $x\in L^{\infty}(\Omega).$

For  $x\in A,$ take a sequence of idempotents  $e_n, n\in\mathbb{N}$
from $L^{\infty}(\Omega)$  such that $e_n x\in L^{\infty}(\Omega)$
and $e_n\uparrow\textbf{1}.$  Then $e_n D(x)=D(e_n x)-D(e_n)x=0,$
i.e. $e_n D(x)=0$ for all $n\in\mathbb{N}.$ Since
$e_n\uparrow\textbf{1}$ this implies that $D(x)=0.$ The proof is
complete. $\blacksquare$

 \vspace{1cm}

\textbf{Acknowledgments.} \emph{The  authors would like to
acknowledge the hospitality of the $\,$ "Institut f\"{u}r Angewandte
Mathematik",$\,$ Universit\"{a}t Bonn (Germany). This work is
supported in part by the DFG 436 USB 113/10/0-1 project (Germany)
and the Fundamental Research Foundation of the Uzbekistan Academy of
Sciences.}

\end{document}